\documentclass[a4paper,reqno,11pt]{amsart}
\usepackage{amsfonts}

\begin{document}
\theoremstyle{plain}
\begingroup
\newtheorem{theorem}{Theorem}[section]
\newtheorem{lemma}[theorem]{Lemma}
\newtheorem{proposition}[theorem]{Proposition}
\newtheorem{corollary}[theorem]{Corollary}
\endgroup

\theoremstyle{definition}
\begingroup
\newtheorem{definition}[theorem]{Definition}
\newtheorem{remark}[theorem]{Remark}
\newtheorem{example}[theorem]{Example}
\endgroup

\theoremstyle{remark}
\begingroup

\endgroup

\mathsurround=1pt
\mathchardef\emptyset="001F

\numberwithin{equation}{section}

\newcommand{\R}{\mathbb R}
\newcommand{\N}{\mathbb N}
\newcommand{\Z}{\mathbb Z}
\newcommand{\leb}{{\mathcal L}}
\newcommand{\mthree}{{\mathbb M}^{3{\times}3}}
\newcommand{\mtwo}{{\mathbb M}^{2{\times}2}}
\newcommand{\dist}{{\rm dist}}
\newcommand{\sym}{{\rm sym}\,}
\newcommand{\skw}{{\rm skew}\,}
\newcommand{\wto}{\rightharpoonup}
\newcommand{\eps}{\varepsilon}
\newcommand{\id}{{\mathit I\!d}}
\newcommand{\ol}{\overline}
\newcommand{\hs}{\hspace{.2cm}}

\title[Convergence of equilibria of planar thin elastic beams]
{Convergence of equilibria of planar\\ thin elastic beams}
%%{\\ \bf DRAFT - NOT FOR CIRCULATION}}
\author{M.G.\ Mora}
\author{S.\ M\"uller}
\author{M.G.\ Schultz}
\address[M.G.~Mora]{SISSA, Via Beirut 2-4, 34014 Trieste, Italy}
\address[S.~M\"uller]{Max Planck Institute for Mathematics in the Sciences, 
Inselstrasse 22, 04103 Leipzig, Germany}
\address[M.G.~Schultz]{Max Planck Institute for Mathematics in the Sciences, 
Inselstrasse 22, 04103 Leipzig, Germany}
\email[Maria Giovanna Mora]{mora@sissa.it}
\email[Stefan M\"uller]{sm@mis.mpg.de}
\email[Maximilian Schultz]{schultz@mis.mpg.de}
%% http://atomoptic.iota.u-psud.fr/

\begin{abstract} We consider a thin elastic strip $\Omega_h = (0,L) \times
(-h/2, h/2)$, and we show that stationary points of the nonlinear elastic
energy (per unit height) $E^h(v) = \frac{1}{h} \int_{\Omega_h} (W(\nabla v) - h^2
g(x_1) {\,\cdot\,} v) \, dx$ whose energy is bounded by $C h^2$ converge to stationary
points of the Euler-Bernoulli functional $J_2(\bar{y}) = \int_0^L (\frac{1}{24}
E \kappa^2 - g{\,\cdot\,} \bar{y}) \, dx_1$ where $\bar{y}: (0,L) \to \R^2$, with
$\bar{y}' = \binom{ \cos \theta}{\sin \theta}$, and where $\kappa = \theta'$.
This corresponds to the equilibrium equation $-\frac{1}{12} E \theta'' +
\tilde{g}{\,\cdot\,} \binom{-\sin \theta}{\cos \theta} = 0$, where $\tilde{g}$ is the
primitive of $g$. The proof uses the rigidity estimate for low-energy
deformations \cite{FJM02} and a compensated compactness argument in a singular
geometry. In addition, possible concentration effects are ruled out by a
careful truncation argument.  \end{abstract}

%\today

\maketitle

{\small

\bigskip

}

\bigskip
\bigskip

\section{Introduction and main result}

The relation between three-dimensional nonlinear elasticity and theories for
lower-dimensional objects such as rods, beams, membranes, plates, and shells
has been an outstanding question since the very beginning of research in
elasticity. In fact, there is a large variety of lower-dimensional theories.
They are usually obtained  by making certain strong a-priori
assumptions on the form of the solutions of the full three-dimensional problem,
and hence their rigorous range of validity is typically unclear.  As
highlighted already in the work of Fritz John, the geometric nonlinearity in
elasticity, i.e., the invariance of the elastic energy under rotations, is one
of the key points. In particular, thin elastic objects can undergo large
rotations even under small loads, and this prevents any analysis based on a
na{\"\i}ve linearization. The first rigorous results were only obtained in the
early 90's using a variational approach that guarantees convergence of
minimizers to a suitable limit problem. In this paper, we discuss the
convergence of possibly non-minimizing stationary points of the elastic energy
functional. 

To set the stage, let us first review the variational setting. 
Consider a cylindrical domain 
$\Omega_h = S{\times}(-\frac{h}{2}, \frac{h}{2})$ where $S$ is a bounded
subset of $\R^2$ with Lipschitz boundary. 
To a deformation $v: \Omega_h \to \mathbb{R}^3$, we associate
the elastic energy (per unit height)
$$
E^h (v) = \frac{1}{h} \int_{\Omega_h} W( \nabla v) \, dz.
$$
We assume that the stored-energy density $W$ 
satisfies the following conditions:
\begin{eqnarray}
& &
W(RF) = W(F)  \quad \forall R \in SO(3) \qquad \mbox{(frame indifference)}, 
\label{eq:H1} 
\\
& & 
W=0 \quad \mbox{on } SO(3),  \label{eq:H2}
\\
& & 
W(F) \geq c\,  \dist^2(F,SO(3)), \quad c>0, \label{eq:H3}
\\
& & 
W \mbox{ is } C^2 \mbox{ in a neighbourhood of } SO(3). 
\label{eq:H4}
\end{eqnarray}
Here, $SO(3)$ denotes the group of proper rotations. 
The frame indifference implies the existence of a function $\tilde{W}$
defined on symmetric matrices such that $W(\nabla v) = \tilde{W}((\nabla v)^T \nabla v)$,
i.e., the elastic energy depends only on the pull-back metric of $v$.

{}For the discussion of the limiting behavior of $E^{h}$ as $h \to 0$, it is convenient to rescale
$\Omega_h$ to a fixed domain $\Omega = S{\times}(-\frac12, \frac12)$
 by a change
of variables, $z = (x_1, x_2, h x_3)$ and $y(x) = v(z)$. With
the notation
$$
\nabla_h y = (\partial_1 y, \partial_2 y, \tfrac{1}{h} \partial_3 y) = (\nabla' y, \tfrac{1}{h}
\partial_3 y),
$$
we thus have 
$$
E^h(v) = I^h(y) = \int_\Omega W(\nabla_h y) \, dx.
$$
The variational approach leads to a hierarchy of limiting theories
depending on the scaling of $I^h$. More precisely, as $h \to 0$ 
$$
\frac{1}{h^\beta} I^h \stackrel{\Gamma}{\longrightarrow} I_\beta
$$
in the sense of $\Gamma$-convergence.
This implies, roughly speaking, 
that minimizers of $I^h$ (subject to suitable
boundary conditions or body forces) converge to minimizers
of $I_\beta$, provided $I^h$ evaluated on the minimizers is bounded by 
$C h^\beta$. $\Gamma$-convergence was first established by
LeDret and Raoult for $\beta = 0$ (see \cite{LR95}), then
for all $\beta \geq 2$ in \cite{FJM02, FJM05} (see also
\cite{Pa01, Pa03} for results for $\beta=2$ under additional
conditions). For $0 < \beta < 5/3$ convergence was recently 
obtained by Conti and Maggi in \cite{CM05}, see also \cite{Co03}.
The exponent $\beta = 5/3$ is conjectured to be relevant
for the crumpling of elastic sheets (see \cite{Witten,Ve04,CM05}). 

Here, we focus on the case $\beta=2$, which leads to Kirchhoff's
geometrically nonlinear bending theory. The natural class $\mathcal{A}$ of admissible functions for the limit problem is given by 
isometric $W^{2,2}$ immersions from $S$ into $\R^3$, i.e., 
$$
\mathcal{A} := \left\{ y \in W^{2,2}(\Omega, \R^3):\ 
\partial_3 y = 0,\  (\nabla' y)^T \nabla' y = Id 
\right\}.
$$
The limiting energy functional is
$$
I_2(y) = 
\left\{ 
\begin{array}{ll}
\displaystyle\frac{1}{24} \int_S  Q_2 (A) \, dx_1 dx_2 & \mbox{if } y \in \mathcal{A}, 
\smallskip 
\\
+ \infty & \mbox{else.}
\end{array}
\right.  
$$
Here, $A$ is the second fundamental form, and $Q_2$ is a quadratic
form that can be computed from the linearization
$D^2 W(Id)$ of the 3d energy at the identity. 
If $W(F) = \frac{1}{2} \dist^2(F, SO(3))$, then simply $Q_2(A) = |A|^2$.

In this paper, we consider the convergence of {\it equilibria} 
for the case $\beta = 2$. 
Instead of treating the full problem of a reduction from 3d to 2d, we focus on the 
simpler case 2d to 1d. Thus, we start from a thin strip
\begin{equation}
\Omega_h  = (0,L){\times}(-\tfrac h2, \tfrac h2),
\end{equation}
and after the rescaling $(z_1, z_2) = (x_1, h x_2)$,
$\nabla_h = (\partial_1, \frac{1}{h} \partial_2)$, we consider the 
functional
$$
J^h(y) = \int_\Omega \big(W(\nabla_h y)  - h^2 g(x_1){\,\cdot\,} y \big) \, dx .
$$
The corresponding $\Gamma$-limit is given by
$$
J_2(\bar y) = \int_0^L \big(\frac{1}{24} E \kappa^2 - g {\,\cdot\,} \bar{y}\big)
 \, dx_1,
$$
where
$$
\bar{y}: (0,L) \to \R^2, \quad \bar{y}' = 
\binom{\cos \theta}{\sin \theta}, \quad
\kappa = \theta'.
$$
The functional $J_2$ takes the value $+ \infty$, if $\bar{y}$ is not of the above
form (here we took the liberty to identify maps  on $\Omega$ which
are independent of $x_2$ with maps on $(0,L)$).
It is convenient to fix one endpoint by requiring $\bar{y}(0) = 0$.
%% and $\bar y'(0)=e_1$,
%%where $\{e_1,e_2\}$ is the canonical basis in $\R^2$.
%%%MGM:added boundary condition on y'
%%%SM: took BC for y' out again and moved
%%%BC for theta into the theorem
%%%SM: Rephrased following sentence. The old one had
%%%SM: a 'dangling particple'
%%
Integrating the linear term by parts, we obtain the following
Euler-Lagrange equation corresponding to the limit functional
%%%is given by 
\begin{equation} \label{eq:EL}
- \frac{1}{12} E \theta^{''} + \tilde{g} {\,\cdot\,} 
\binom{- \sin \theta}{ \cos \theta} = 0,
\quad
\tilde{g}(x_1) := \int_L^{x_1} g(\xi) \, d\xi. 
\end{equation}
%%%
%%\begin{equation} \label{eq:ELbc}
%%\theta(0) = 0, \quad \theta'(L) = 0.
%%%\end{equation}

%%%MGM:added this sentence
The main result of the paper is the following.

\begin{theorem}\label{thm:1.1}
Assume (\ref{eq:H1})--(\ref{eq:H4}), the energy $W$ is differentiable,
and the derivative $DW$ is globally Lipschitz. Let $(y^{(h)})$ be a sequence
of stationary points of $J^h$, subject to the boundary condition
$y^{(h)}(0,x_2) = (0,hx_2)$ at $x_1=0$ and to
natural boundary conditions on the remaining
boundaries. Assume further 
\begin{equation} \label{eq:energy}
\int_\Omega W(\nabla_h y^{(h)}) \leq C h^2.
\end{equation}
Then, up to subsequences,
\begin{equation}\label{th1}
y^{(h)} \to \bar{y} \quad \mbox{in } W^{1,2}(\Omega; \mathbb{R}^2),
\end{equation}
as $h\to 0$. The limit function $\bar{y}$ satisfies
\begin{equation}\label{th2}
\partial_1 \bar{y} = 
\binom{\cos \theta}{ \sin \theta}, \quad \partial_2 \bar{y} = 0, 
\end{equation}
and $\theta$ satisfies  (\ref{eq:EL}) and 
%%%%(\ref{eq:ELbc})
\begin{equation} \label{eq:ELbc}
\theta(0) = 0, \quad \theta'(L) = 0.
\end{equation}
\end{theorem}

\begin{remark}
An easy application of Poincar\'e's inequality shows that the estimate
(\ref{eq:energy}) automatically holds for minimizers.
\end{remark}

\begin{remark}
In \cite{Mi88}, Mielke uses a center manifold approach to
compare solutions in a thin strip to a 1d problem. 
His approach already works for finite $h$,
but requires that the nonlinear strain $(\nabla_h y)^T \nabla_h y$
is  $L^\infty$-close to the identity. Applied forces $g$
are also difficult to handle. 
\end{remark}

%%%SM
One key idea in 
the proof of Theorem~\ref{thm:1.1} 
is to replace the use of comparison functions in $\Gamma$-convergence
by a compensated compactness argument (see Steps~6 and 7 in
the next section). To set the stage for this argument, we
use the quantitative rigidity estimate of
\cite{FJM02} to introduce suitable strain-like and stress-like 
variables $G^{(h)}$ and $E^{(h)}$, which are almost curl-free
and divergence-free, respectively (see Steps~2 and 3 below).
To control possible concentration effects, we
use a truncation argument, see Section 3.

\section{Proof}
\label{sec:Proof}
Let $(y^{(h)})$ be a sequence of stationary points of $J^h$,
 i.e., suppose that
\begin{equation}\label{EL1}
\int_\Omega \Big( DW(\nabla_h y^{(h)}){\,:\,}\nabla_h\psi -h^2g{\,\cdot\,}\psi\Big)\, dx =0
\end{equation}
for every $\psi\in W^{1,2}(\Omega;\R^2)$ with $\psi=0$ on $\{x_1=0\}$. 
Assume further estimate (\ref{eq:energy}).

\medskip

\noindent
{\em Step 1. Decomposition of the
deformation gradients in rotation and strain.}
\smallskip

\noindent
By Proposition~\ref{rotation}, we can construct a sequence $(R^{(h)})\subset C^\infty((0,L);\mtwo)$
such that for every $x_1 \in (0,L)$, $R^{(h)}(x_1)\in SO(2)$ and
\begin{eqnarray}
&
\|\nabla_h y^{(h)}-R^{(h)}\|_{L^2}\le Ch, \label{rig1}
\\
&
\|(R^{(h)})'\|_{L^2}+h \|(R^{(h)})''\|_{L^2}\le C, \label{rig2} \\
& |R^{(h)}(0) - Id | \leq C \sqrt{h}. \label{rig3}
\end{eqnarray}
By (\ref{rig2}), there exists $R\in W^{1,2}((0,L);\mtwo)$ such that up to subsequences $R^{(h)}$
converges to $R$ weakly in $W^{1,2}((0,L);\mtwo)$, hence uniformly in $L^\infty((0,L);\mtwo)$. 
Thus $R(x_1)\in SO(2)$ for every $x_1\in(0,L)$.
Moreover, estimate (\ref{rig1}) implies that
$$
\nabla_h y^{(h)}\to R \quad \mbox{strongly in } L^2(\Omega;\mtwo).
$$
In particular, $\partial_2 y^{(h)} \to 0$ and thus
\begin{equation}\label{gradient}
\nabla y^{(h)}\to Re_1\otimes e_1  \quad \mbox{strongly in } L^2(\Omega;\mtwo).
\end{equation}
Since $|y^{(h)}(0,x_2)| \leq h \to 0$, we deduce from
Poincar\'e's inequality that $y^{(h)} \to \bar{y}$
strongly in $W^{1,2}(\Omega;\R^2)$ and that 
$\bar{y}$ satisfies $\partial_1\bar{y}=Re_1$, 
$\partial_2\bar{y}=0$ a.e.\ in $\Omega$.
Thus, (\ref{th1}) and (\ref{th2}) are proved.

We now make use of the approximated sequence of rotations $R^{(h)}$ in order to decompose the deformation
gradients as 
\begin{equation}\label{decomp}
\nabla_h y^{(h)}=R^{(h)}(Id+hG^{(h)}).
\end{equation}
The 
 $G^{(h)}:\Omega\to\mtwo$ are bounded in $L^2(\Omega;\mtwo)$ by (\ref{rig1}).
Thus, up to extracting a subsequence, we can assume 
\begin{equation}\label{Gdeb}
G^{(h)}\wto G \quad \mbox{weakly in } L^2(\Omega;\mtwo)
\end{equation}
for some $G\in L^2(\Omega;\mtwo)$. 

\medskip

\noindent
{\em Step 2. Consequences of the compatibility of the strain.}
\smallskip

\noindent
Up to the factor $(R^{(h)})^T$, the strains $G^{(h)}$ are essentially scaled gradients.
This has some important consequences on the form of the limit strain $G$.
To deduce these properties, it is convenient to introduce a sequence of auxiliary functions 
$z^{(h)}:\Omega\to\R^2$ defined as
\begin{equation}\label{def:z^h}
z^{(h)}(x):=\frac1h  y^{(h)}(x) -\frac1h \int_0^{x_1} R^{(h)}(s)e_1\, ds
- x_2R^{(h)}(x_1)e_2 .
\end{equation}
By the definition (\ref{decomp}) of $G^{(h)}$,
\begin{eqnarray}
\nabla_h z^{(h)} & = & \frac1h\Big(\nabla_h y^{(h)} -R^{(h)}\Big) - x_2(R^{(h)})'e_2\otimes e_1
\nonumber
\\
& = & R^{(h)}\Big(G^{(h)}- x_2(R^{(h)})^T(R^{(h)})'e_2\otimes e_1  \Big).
\label{nablaz}
\end{eqnarray}
Since $R^{(h)} \in SO(2)$, there exist
 $\theta^{(h)}\in C^\infty(0,L)$ such that
\begin{equation} \nonumber
R^{(h)} = \left( 
  \begin{array}{cc}  \cos\theta^{(h)} &  -\sin\theta^{(h)} \\
                     \sin\theta^{(h)} & \cos\theta^{(h)}
  \end{array}
          \right) \, .
\end{equation}
Then $(R^{(h)})^T(R^{(h)})'e_2=- (\theta^{(h)})'e_1$, and
 equality (\ref{nablaz}) can be rewritten as
\begin{equation}\label{form1}
\nabla_h z^{(h)} = R^{(h)}\Big(G^{(h)}+ x_2(\theta^{(h)})'e_1\otimes e_1  \Big).
\end{equation}
Now, recall that $R^{(h)}$ converges uniformly to $R$, 
$G^{(h)}$ converges weakly to $G$ in $L^2(\Omega;\mthree)$,  
and $(\theta^{(h)})'$ converges weakly to $\theta'$ in $L^2(0,L)$, 
where $\theta$ satisfies (\ref{th2}).
{}From these properties it follows that
\begin{equation}\label{nablaz2}
\nabla_h z^{(h)} \wto R\Big(G+ x_2\theta'e_1\otimes e_1\Big)
\quad \mbox{weakly in } L^2(\Omega;\mtwo).
\end{equation}
The definition of $z^{(h)}$ and (\ref{rig3}) yield 
$|z^{(h)}(0,x_2)| \leq C \sqrt{h}$. Hence Poincar\'e's inequality
shows that $z^{(h)}$ converges weakly in  $W^{1,2}(\Omega;\R^2)$.
The limit function $z$ satisfies
\begin{equation}\label{form2}
R^T\partial_1 z= Ge_1+x_2\theta'e_1, \quad \partial_2 z=0 \quad \hbox{a.e.\ in } \Omega.
\end{equation}
In particular, $z$ does not depend on $x_2$, and thus, 
by the first equality in (\ref{form2}), the vector 
$Ge_1$ is linear in $x_2$.

Let $\hat G:(0,L)\to\mtwo$ be the first moment of $G$ defined by 
$$
\hat G(x_1):=\int_{-\frac12}^{\frac12}x_2G(x)\, dx_2.
$$ 
As $R$ and $z$ are independent of $x_2$, we deduce from (\ref{form2}) that
\begin{equation}\label{hatG}
\hat G_{11}=-\frac{1}{12}\theta', \quad \hat G_{21}=0 \quad \hbox{a.e.\ in } (0,L).
\end{equation}

\medskip

\noindent
{\em Step 3. Consequences of the Euler-Lagrange equations.}
\smallskip

\noindent
Let $E^{(h)}:\Omega\to\mtwo$ be the scaled stress defined by
\begin{equation}\label{Ehdef}
E^{(h)}:=\frac1h DW(Id+h G^{(h)}).
\end{equation}
Since $DW$ is Lipschitz continuous and the
$G^{(h)}$ are bounded in $L^2(\Omega;\mtwo)$, 
the functions $E^{(h)}$ are bounded in $L^2(\Omega;\mtwo)$.
By Proposition~\ref{Ehweak}, we have 
\begin{equation}\label{lebG}
E^{(h)}\wto E:=\leb G \quad \mbox{weakly in }L^2(\Omega;\mtwo),
\end{equation}
where the linear map $\leb$ on the matrix space is given by
$\leb:=D^2W(Id)$. We note in particular that $E$ is symmetric
since, by frame indifference, $\leb F = \leb \sym F$ and 
$\leb F = (\leb F)^T$ for all $F \in \mtwo$.

By the decomposition (\ref{decomp}) and frame indifference of $W$, we obtain 
$$
DW(\nabla_h y^{(h)})= R^{(h)} DW(Id+hG^{(h)}) = h R^{(h)}E^{(h)}.
$$
The Euler-Lagrange equations (\ref{EL1}) can be written in terms of the stresses $E^{(h)}$:
\begin{equation}\label{EL2}
\int_\Omega \Big( R^{(h)}E^{(h)}{\,:\,}\nabla_h\psi -hg{\,\cdot\,}\psi\Big)\, dx =0
\end{equation}
for every $\psi\in W^{1,2}(\Omega;\R^2)$ with $\psi=0$ on $\{x_1=0\}$. 
Multiplying (\ref{EL2}) by $h$ and passing to the limit as $h\to 0$, we find
\begin{equation}\label{cons1}
\int_\Omega REe_2{\,\cdot\,}\partial_2\psi \, dx =0
\end{equation}
for every $\psi\in W^{1,2}(\Omega;\R^2)$ with $\psi=0$ on $\{x_1=0\}$. 
This yields $REe_2=0$ a.e.\ in $\Omega$ and hence 
 $Ee_2=0$ a.e.\ in $\Omega$. Therefore,
as $E$ is symmetric, we conclude that
\begin{equation}\label{cons2}
E(x)= E_{11}(x)e_1\otimes e_1
\end{equation}
for a.e.\ $x\in\Omega$.

\medskip

\noindent
{\em Step 4. Symmetry properties of $E^{(h)}$.}
\smallskip

\noindent
Since $W$ is frame indifferent, the matrix $DW(F)F^T$ is symmetric.
Choosing $F=Id+hG^{(h)}$, we deduce that 
\begin{equation}\label{symm}
E^{(h)}-(E^{(h)})^T= -h \big( E^{(h)}(G^{(h)})^T-G^{(h)}(E^{(h)})^T \big).
\end{equation}
Using the boundedness of $E^{(h)}$ and $G^{(h)}$ in $L^2(\Omega;\mtwo)$,
we have in particular the estimate
\begin{equation}\label{cons3}
\|E^{(h)}_{12}- E^{(h)}_{21}\|_{L^1}\le Ch. 
\end{equation}

\medskip

\noindent
{\em Step 5. Moments of the Euler-Lagrange equations.}
\smallskip

\noindent
Let us introduce the zeroth and the first moment of the stress $E^{(h)}$,
$$
\bar E^{(h)}(x_1):=\int_{-\frac12}^{\frac12}E^{(h)}(x)\,dx_2, \qquad
\hat E^{(h)}(x_1):=\int_{-\frac12}^{\frac12} x_2E^{(h)}(x)\,dx_2,
$$
for every $x_1\in (0,L)$. In the following, we will derive the equations satisfied 
by these moments.

Let $\varphi\in C^\infty([0,L];\R^2)$ be such that $\varphi(0)=0$. Using $\varphi$ as
a test function for the Euler-Lagrange equation (\ref{EL2}), we obtain
$$
\int_\Omega \Big( R^{(h)}E^{(h)}e_1{\,\cdot\,}\varphi' -hg{\,\cdot\,}\varphi\Big)\, dx =0.
$$
As $R^{(h)}$, $\varphi$, and $g$ depend only on the variable $x_1$, this 
equality is equivalent to
$$
\int_0^L \Big( R^{(h)}\bar E^{(h)}e_1{\,\cdot\,}\varphi' -hg{\,\cdot\,}\varphi\Big)\, dx_1 =0.
$$
This equation holds for every $\varphi\in C^\infty([0,L];\R^2)$ with $\varphi(0)=0$,
and hence we deduce that 
\begin{equation}\label{mom0}
\bar E^{(h)}e_1= -h(R^{(h)})^T \tilde g \quad \hbox{a.e.\ in } (0,L),
\end{equation}
where $\tilde g$ is the primitive of $g$ defined in (\ref{eq:EL}). By passing
to the limit, we obtain
\begin{equation}\label{mom1}
\bar Ee_1=0 \quad \hbox{a.e.\ in } (0,L).
\end{equation}

As for the first moment, 
let $\varphi\in C^\infty([0,L];\R^2)$ be again such that $\varphi(0)=0$.
Using $\psi(x):=x_2\varphi(x_1)$ as a test function in (\ref{EL2}), we now obtain
$$
\int_\Omega \Big( x_2R^{(h)}E^{(h)}e_1{\,\cdot\,}\varphi'+\frac1h R^{(h)}E^{(h)}e_2{\,\cdot\,}\varphi -hx_2
g{\,\cdot\,}\varphi\Big)\, dx =0.
$$
Upon integration with respect to $x_2$, this equation becomes
$$
\int_0^L \Big(R^{(h)} \hat E^{(h)}e_1{\,\cdot\,}\varphi' 
+\frac1h R^{(h)} \bar E^{(h)}e_2{\,\cdot\,}\varphi\Big)\, dx_1 =0.
$$
We can choose $\varphi$ to be of the form $\varphi=\phi\,R^{(h)} e_1$
with $\phi\in C^\infty([0,L])$ and $\phi(0)=0$. Hence
$$
\int_0^L \Big( \hat E^{(h)}_{11}\phi' + \hat E^{(h)}e_1{\,\cdot\,} \phi (R^{(h)})^T(R^{(h)})'e_1
+\frac1h  \bar E^{(h)}_{12}\phi\Big)\, dx_1 =0.
$$
Wtih  the identity $(R^{(h)})^T(R^{(h)})'e_1= (\theta^{(h)})'e_2$, 
this expression reduces to
\begin{equation}\label{EL5}
\int_0^L \Big( \hat E^{(h)}_{11}\phi' +(\theta^{(h)})'\hat E^{(h)}_{21}\phi
+\frac1h  \bar E^{(h)}_{12}\phi\Big)\, dx_1 =0.
\end{equation}

{}From the estimate (\ref{cons3}) and the identity (\ref{mom0}), we infer an $L^1(0,L)$-bound on $\frac{1}{h}\bar E_{12}^{(h)}$.
Since $(\theta^{(h)})'$ and $\hat E^{(h)}_{21}$ are bounded in $L^2(0,L)$, 
$(\theta^{(h)})'\hat E^{(h)}_{21}$ is clearly bounded in $L^1(0,L)$.
Therefore, equation (\ref{EL5}) implies that 
\begin{equation}\label{EL6}
\|\partial_1\hat E^{(h)}_{11}\|_{L^1}\le C, \quad \hat E^{(h)}_{11}(L)=0.
\end{equation}
Hence, the sequence $\hat E^{(h)}_{11}$ is strongly compact in $L^p(0,L)$ for every $p<\infty$.

\medskip

\noindent
{\em Step 6. Convergence of the energy by the div-curl lemma.}
\smallskip

\noindent
The strong compactness of the sequence $(\hat E^{(h)}_{11})$ allows us to pass to the limit
in the energy integral
$$
\frac{1}{h^2}\int_\Omega DW(Id+hG^{(h)}){\,:\,}hG^{(h)}\, dx = 
\int_\Omega E^{(h)}{\,:\,}G^{(h)}\, dx.
$$
The limit is obtained by exploiting the div-curl structure of the product $E^{(h)}{\,:\,}G^{(h)}$. 
By the Euler-Lagrange equation (\ref{EL2}), the scaled divergence of $R^{(h)}E^{(h)}$ is infinitesimal
in $L^2(\Omega;\R^2)$ as $h\to 0$, while by the decomposition (\ref{form1}), the matrix
$R^{(h)}G^{(h)}$ has essentially the structure of a scaled gradient.

Let $\varphi\in C^\infty([0,L])$ with $\varphi(0) = 0$. 
Using formula (\ref{form1}), we have
\begin{eqnarray}
\lefteqn{\int_\Omega \varphi E^{(h)}{\,:\,}G^{(h)}\, dx = 
\int_\Omega \varphi R^{(h)}E^{(h)}{\,:\,}R^{(h)}G^{(h)}\, dx} \nonumber
\\
& = & \int_\Omega \varphi R^{(h)}E^{(h)}{\,:\,}\nabla_h z^{(h)}\, dx
-\int_\Omega \varphi E^{(h)}_{11}x_2(\theta^{(h)})'\, dx. \label{div1}
\end{eqnarray}
We deduce from the Euler-Lagrange equation (\ref{EL2}), applied
with $\psi = \varphi z^{(h)}$, that
$$
\int_\Omega \varphi R^{(h)}E^{(h)}{\,:\,}\nabla_h z^{(h)}\, dx =
h\int_\Omega \varphi g{\,\cdot\,} z^{(h)}\, dx 
- \int_\Omega \varphi'R^{(h)}E^{(h)}e_1{\,\cdot\,} z^{(h)}\, dx.
$$
Since $z^{(h)}$ converges strongly to $z$ in $L^2(\Omega; \R^2)$, we can pass to the limit
in the above formula and obtain
\begin{eqnarray}
\lim_{h\to 0}\int_\Omega \varphi R^{(h)}E^{(h)}{\,:\,}\nabla_h z^{(h)}\, dx & = &
- \int_\Omega \varphi'RE e_1{\,\cdot\,} z\, dx \nonumber
\\
& = &
- \int_0^L \varphi'R\bar E e_1{\,\cdot\,} z\, dx_1 \ = \ 0, \label{div2}
\end{eqnarray}
where the last two equalities follow from the independence of $z$ of $x_2$ and the identity (\ref{mom1}).

As for the last term in (\ref{div1}), we have
$$
\int_\Omega \varphi E^{(h)}_{11}x_2(\theta^{(h)})'\, dx=
\int_0^L \varphi \hat E^{(h)}_{11} (\theta^{(h)})'\, dx_1
$$
and, using the strong convergence of $\hat E^{(h)}_{11}$ proved at the end of Step~5,
we deduce that
\begin{equation}\label{div3}
\lim_{h\to 0}\int_\Omega \varphi E^{(h)}_{11}x_2(\theta^{(h)})'\, dx=
\int_0^L \varphi \hat E_{11} \theta'\, dx_1.
\end{equation}
{}From the first equality in (\ref{form2}), we infer that
$$
\int_0^L \varphi \hat E_{11} \theta'\, dx_1 =
\int_\Omega \varphi E_{11} x_2\theta'\, dx =  \int_\Omega \varphi E_{11} R^Tz'{\,\cdot\,} e_1\, dx
- \int_\Omega \varphi E_{11}G_{11}\, dx.
$$
Using (\ref{mom1}) and (\ref{cons2}), the previous equality yields
\begin{eqnarray}
\int_0^L \varphi \hat E_{11} \theta'\, dx_1 & = & \int_0^L \varphi \bar E_{11} R^Tz'{\,\cdot\,} e_1\, dx_1
-\int_\Omega \varphi E_{11}G_{11}\, dx
\nonumber
\\
& = & -\int_\Omega \varphi E{\,:\,}G \, dx.
\label{div4}
\end{eqnarray}
{}Finally combining equations (\ref{div1})--(\ref{div4}), 
we obtain convergence of the energies
\begin{equation}\label{en-conv}
\lim_{h\to 0}\int_\Omega \varphi E^{(h)}{\,:\,}G^{(h)}\, dx
= \int_\Omega \varphi E {\,:\,}G \, dx
\end{equation}
for every $\varphi\in C^\infty([0,L])$ with $\varphi(0) = 0$.

\medskip

\noindent
{\em Step 7. Strong convergence of the symmetric part of $G^{(h)}$.}
\smallskip

\noindent
In order to emphasize the structure of the argument, we first conclude
the proof under the additional assumption 
\begin{equation}\label{hGh}
\lim_{h\to 0}\, h \| G^{(h)}\|_{L^\infty}=0.
\end{equation}
This assumption will allow us to replace $E^{(h)}$ in (\ref{en-conv})
by $\leb G^{(h)}$. Since we already know that $E = \leb G$
and since $\leb$ is positive definite on symmetric matrices,
we can conclude strong convergence of $G^{(h)}$ and hence of
$E^{(h)}$ (away from $x_1 = 0$). Using this strong convergence,
we can easily pass to the limit in (\ref{EL5}), using 
formula (\ref{symm}) for the skew-symmetric part of $E^{(h)}$ and the
good control on $E^{(h)}_{21}$. 
We will show in the next section how assumption (\ref{hGh})
can be avoided
through the use of a careful truncation argument.

A Taylor expansion of $DW$ around the identity matrix yields
\begin{equation}\label{Eexp}
E^{(h)}= \frac1h DW(Id+hG^{(h)})= \leb G^{(h)} +\frac1h\eta(hG^{(h)}),
\end{equation}
where $|\eta(A)|/|A|\to 0$, as $|A|\to 0$. For every $t>0$, we define
$$
\omega(t):=\sup\Big\{\frac{|\eta(A)|}{|A|}: \ |A|\le t \Big\}.
$$
Then  $\omega(t)\to 0$, as $t\to 0^+$.
Using the expansion (\ref{Eexp}), we obtain
$$
|E^{(h)}{\,:\,}G^{(h)}-\leb G^{(h)}{\,:\,}G^{(h)}|\le
\frac1h|\eta(hG^{(h)})|\,|G^{(h)}|\le
\omega(h \|G^{(h)}\|_{L^\infty})|G^{(h)}|^2.
$$
Assumption (\ref{hGh}) implies that the last term in this inequality
tends to zero in $L^1(\Omega)$, as $h\to 0$. By (\ref{en-conv}) and (\ref{lebG}),
we thus obtain for every $\varphi\in C^\infty([0,L])$
with $\varphi(0) = 0$
\begin{equation}\label{101}
\lim_{h\to 0}\int_\Omega \varphi \leb G^{(h)}{\,:\,}G^{(h)}\, dx
= \int_\Omega \varphi \leb G {\,:\,}G \, dx.
\end{equation}

{}From the assumptions on $W$, we obtain a constant $C>0$ such that
$$
\leb A{\,:\,}A\geq C|\sym A|^2
$$ 
for every $A\in\mtwo$. This inequality, together with (\ref{Gdeb}) and (\ref{101}), implies that
$$
\lim_{h\to 0} \, C\! \int_\Omega \varphi \, |\sym (G^{(h)}-G) |^2 \, dx \leq
\lim_{h\to 0} \int_\Omega \varphi \, \leb (G^{(h)}-G){\,:\,}(G^{(h)}-G) \, dx =0
$$
for every nonnegative $\varphi\in C^\infty([0,L])$ with $\varphi(0) = 0$.
In particular, we have
\begin{equation}\label{en-conv1}
\sym G^{(h)}\to \sym G \quad \mbox{strongly in } 
L^2((a,L) {\times}(-\tfrac12, \tfrac12);\mtwo)
\end{equation}
for every $a>0$.
Since $\leb A=\leb\sym A$ for every $A\in\mtwo$, the Taylor expansion (\ref{Eexp}),
together with (\ref{hGh}) and (\ref{en-conv1}), yields
\begin{equation}\label{en-conv2}
E^{(h)}\to E=\leb\sym G \quad \mbox{strongly in } 
L^2((a,L) {\times}(-\tfrac{1}{2},\tfrac{1}{2});\mtwo)
\end{equation}
for every $a>0$.

\medskip

\noindent
{\em Step 8. Derivation of the limit equation.}
\smallskip

\noindent
Due to the convergence (\ref{en-conv2}), we can pass to the limit in (\ref{symm})
and obtain
\begin{equation}\label{Eskew}
\frac1h \big(E^{(h)}-(E^{(h)})^T\big)\wto GE^T-EG^T
\end{equation}
weakly in $L^1((a,L){\times}(-\frac12,\frac12);\mtwo)$ for every 
$a>0$.
Note that by (\ref{cons2}) one has
$$
GE^T-EG^T= E_{11}(Ge_1\otimes e_1-e_1\otimes Ge_1)=
E_{11}G_{21}(e_2\otimes e_1-e_1\otimes e_2)
$$
and recall that $G_{21}$ is independent of $x_2$ by (\ref{form2}).
Thus, by (\ref{Eskew}) and (\ref{mom1}) we deduce that
%%%MGM:the convergence is in L^1(a,L)
\begin{equation}\label{Eskew12}
\frac1h ( \bar{E}^{(h)}_{21}- \bar{E}^{(h)}_{12} )\wto G_{21}\bar E_{11}=0
\quad \mbox{weakly in } 
L^1(a,L)
\end{equation}
for every $a>0$.

We now have all the necessary ingredients to derive the limit equation. 
It follows from (\ref{mom0}) that $\frac1h \bar E^{(h)}_{21}$ 
converges to 
$-R^T\tilde g{\,\cdot\,} e_2$ strongly in $L^2(0,L)$. 
Combining this with (\ref{Eskew12}), we deduce that
$\frac1h \bar E^{(h)}_{12}$ converges to $-R^T\tilde g{\,\cdot\,} e_2$ weakly in $L^1(a,L)$ for every $a>0$. 
Using the strong convergence (\ref{en-conv2}) and the fact that $\hat E_{21}=0$ by (\ref{cons2}),
we can pass to the limit in the equation (\ref{EL5}) and conclude that
\begin{equation}\label{ELfin}
\int_0^L \Big( \hat E_{11} \varphi' 
  -  R^T\tilde g{\,\cdot\,} e_2 \varphi\Big)\, dx_1 =0
\end{equation}
for every $\varphi\in C^\infty([0,L])$ which vanishes on an interval
$(0,a)$, with $a>0$. By approximation we see that the
limiting equation holds for all $\varphi \in C^\infty([0,L])$
such that $\varphi(0)=0$.

Since $E=\leb\sym G$ and $E =E_{11}e_1\otimes e_1$, we have
$$
\sym G=E_{11}\leb^{-1}(e_1\otimes e_1),
$$
which yields
$$
G_{11}=E_{11}\leb^{-1}(e_1\otimes e_1){\,:\,}(e_1\otimes e_1).
$$
Using the representation formula (\ref{hatG}), we deduce that 
$$
-\frac{1}{12}\theta'= \hat E_{11}\leb^{-1}(e_1\otimes e_1){\,:\,}(e_1\otimes e_1).
$$
Thus, combining  this equality and (\ref{ELfin}), we conclude that
(\ref{eq:EL}) is satisfied with $E^{-1}=\leb^{-1}(e_1\otimes e_1){\,:\,}(e_1\otimes e_1)$. Moreover, 
the natural boundary condition
follows directly from (\ref{ELfin}),
while (\ref{rig3}) and the uniform convergence of $R^{(h)}$ imply
$\theta(0) =~0$. 
\qed

\section{Truncation and compactness}
\label{sec:trcp}

Apart from some minor issues which are discussed in the next section,
the main point is to remove the strong hypothesis $h\|G^{(h)}\|_\infty\to 0$ in the
argument of Section~\ref{sec:Proof}. To achieve this, we will use a truncation argument.
We first observe that the standard truncation result can also be applied to functions 
defined in thin rectangles, equivalently it can be applied to the scaled gradient 
$\nabla_h=(\partial_1,\frac1h\partial_2)$. Moreover by a good choice of 
the truncation parameter, the bad set where the truncation does not agree with the
original function can be chosen to be particularly small
(see Lemma~\ref{goodtrunc} below).

\begin{lemma} \label{truncate} 
There exists a constant $C$ with the following property.
{}For every $h>0$, every $A > a > 0$ and every $u\in W^{1,2}(\Omega_h;\R^2)$
there exist $\lambda \in [a,A]$ and a function
 $v\in W^{1,\infty}(\Omega_h;\R^2)$ such that
\begin{equation} \label{trunc1}
\| \nabla v \|_{L^\infty} \leq \lambda,
\end{equation}
\begin{equation} \label{trunc2}
\lambda^2 \leb^2(\{x\in \Omega_h: \ u(x)\neq v(x)\})\le
\frac{C}{\ln(A/a)} 
\int_{\Omega_h} |\nabla u|^2\, dx. 
\end{equation}
\end{lemma}

{}For the proof, we refer to Section~\ref{Sec:aux}. We also recall that
\begin{equation}
\nabla u = \nabla v \quad \mbox{a.e.\ in the set }  
\{ x \in \Omega_h : u(x) = v(x) \}.
\end{equation}

\begin{proof}[Conclusion of the proof of Theorem~\ref{thm:1.1}]
Using a truncation, we first define functions $\tilde{y}^{(h)}$
such that the corresponding rescaled strains $\tilde{G}^{(h)}$
satisfy
\begin{equation} \label{htGh}
\lim_{h \to 0} h ||\tilde{G}^{(h)}||_{L^\infty} = 0.
\end{equation}
We can then use a Taylor expansion as in Step 7 of the previous
section to conclude that
$\int \varphi \leb \tilde{G}^{(h)}{\,:\,} \tilde{G}^{(h)} \approx  \int \varphi \tilde{E}^{(h)}{\,:\,}\tilde{G}^{(h)}$.
The crucial step is to show strong convergence of $\sym \tilde{G}^{(h)}$
(away from $x_1 = 0$).
It will be easy to see that $\tilde{G}^{(h)}$ and $G^{(h)}$ have the same
weak limit. The main point is to get strong $L^2$ convergence
of the truncated sequence $\tilde{G}^{(h)}$.
We can adapt the compactness argument in
Step 7 of the previous section to get convergence of the truncated sequence
if we can show that
$$
\int_\Omega \varphi E^{(h)}{\,:\,} G^{(h)} \, dx -
\int_\Omega \varphi \tilde{E}^{(h)} {\,:\,} \tilde{G}^{(h)} \, dx \to 0.
$$
To prove this, we exploit that the
most dangerous term in this difference, namely 
$E^{(h)}{\,:\,}(\tilde{G}^{(h)} - G^{(h)})$ has an (approximate) 
div-curl structure.
{}Finally, we need to pass to the limit in (\ref{EL5}).
The difficulty is that at this point we only know  strong convergence of
$\tilde{E}^{(h)}$ and not of $E^{(h)}$. To control the remainder
term, we first use the fact that $(\theta^{(h)})'$ cannot concentrate 
a finite amount of $L^2$-norm  on the set where the truncation  
deviates from the original function. To estimate the skew-symmetric
part of $E^{h}$ we use its representation (\ref{symm}) in connection with
another application of the div-curl lemma.

\medskip

\noindent{\em Step 1. Definition of the truncated functions.}
\smallskip

\noindent
We consider the functions $z^{(h)}$ defined in (\ref{def:z^h}) and
their rescalings $\check{z}^{(h)}(x) := z^{(h)}(x_1, x_2/h)$.
Applying Lemma~\ref{truncate} to $\check{z}^{(h)}$ with 
$a = h^{-5/8}$, $A= h^{-7/8}$ and undoing the rescaling,
we obtain functions
$\tilde z^{(h)}:\Omega\to\R^2$ and $\lambda_h \in [h^{-5/8}, h^{-7/8}]$
with the following properties:
\begin{equation}\label{tilde1}
 \|\nabla_h \tilde z^{(h)}\|_{L^\infty}\le \lambda_h, 
\end{equation}
\begin{eqnarray}
\lambda_h^2 \leb^2(A_h) & \leq &
\frac{C}{\ln (1/h)}
\int_{\Omega} |\nabla_h z^{(h)}|^2\, dx 
 \nonumber
\\
&  \leq & 
\frac{C}{\ln (1/h)}
\int_\Omega \Big(|G^{(h)}|^2+|(\theta^{(h)})'|^2\Big)\, dx, \label{tilde2}
\end{eqnarray}
where 
$A_h :=\{x\in \Omega: \ z^{(h)}(x)\neq \tilde z^{(h)}(x)\}$.
In particular, we have
\begin{equation}\label{tilde4}
h^{1/2} \lambda_h \to \infty, \quad
h \lambda_h \to 0, \quad \mbox{and} \quad 
\lambda_h^2  \leb^2(A_h)  \to 0.
\end{equation}
We can also define a sequence of approximated deformations $\tilde y^{(h)}:\Omega\to\R^2$ which are associated with the auxiliary functions $\tilde z^{(h)}$:
$$
\tilde y^{(h)}:=h\tilde z^{(h)}+\int_0^{x_1}R^{(h)}(s)e_1\, ds+hx_2R^{(h)}e_2.
$$
Let $\tilde G^{(h)}:\Omega\to\mtwo$ be the corresponding approximated strains defined by
the relation
$$
\nabla_h\tilde y^{(h)} = R^{(h)} (Id + h  \tilde G^{(h)}),
$$
and let $\tilde E^{(h)}:\Omega\to\mtwo$ be the corresponding stresses defined as
\begin{equation}\label{defEtilde}
\tilde E^{(h)}:= \frac1h DW(Id + h  \tilde G^{(h)}).
\end{equation}
Using the definition of $\tilde y^{(h)}$, it is easy to see that
\begin{equation}\label{defGtilde}
\tilde G^{(h)} = (R^{(h)})^T\nabla_h\tilde z^{(h)}-x_2(\theta^{(h)})'e_1\otimes e_1.
\end{equation}
{}From (\ref{tilde1}) and (\ref{tilde2}),
%%%MGM:changed references
we easily see that 
$\nabla_h \tilde{z}^{(h)}$ and $\nabla_h z^{(h)}$ are bounded in 
$L^2$. In fact, they have the same weak limit. To see this, fix
$\eta \in L^2(\Omega;\mtwo)$. Then 
$$ 
\left| \int_\Omega \eta{\,:\,} (\nabla_h  \tilde{z}^{(h)} - \nabla_h z^{(h)}) \, dx
\right|
\le 
C \Big( \int_{A_h}  |\eta|^2 \, dx \Big)^{1/2} \to 0.
$$
Thus 
%%%MGM:added a \tilde on G^h
\begin{equation} \label{tildeweak}
\tilde G^{(h)} \wto G  \quad \mbox{weakly in } L^2(\Omega;\mtwo).
\end{equation}

\medskip

\noindent {\em Step 2.  $L^\infty$-convergence of $h \tilde{G}^{(h)}$
and strong convergence of $\sym \tilde{G}^{(h)}$ and $\tilde{E}^{(h)}$.}
\smallskip

\noindent
We recall the estimate
\begin{equation} \label{interpol}
\sup | f - \bar{f}|^2 \leq 2 ||f||_{L^2} ||f'||_{L^2}, \quad
\mbox{where } \bar{f} = \frac{1}{L} \int_0^L f \, dx, 
\end{equation}
which follows from the identity $(g^2)' = 2 g g'$, applied with $g= f - \bar{f}$.
Since $h(\theta^{(h)})''$ and
$(\theta^{(h)})'$ are bounded in $L^2(0,L)$ by (\ref{rig2}), we 
have in particular that $(1/L) \int_0^L |(\theta^{(h)})'| \leq C$ and thus
(\ref{interpol}) yields  $|(\theta^{(h)})'|\le Ch^{-1/2}$. 
This estimate and inequality (\ref{tilde1}) imply 
\begin{equation}\label{tildeG}
h|\tilde G^{(h)}|\le Ch\lambda_h +Ch^{1/2} \to 0.
\end{equation}
Here, we used that $h\lambda_h\to 0$.
Expanding $DW$ around the identity as in (\ref{Eexp}), we obtain
$$
 \tilde E^{(h)} {\,:\,} \tilde G^{(h)}
= \leb \tilde G^{(h)}{\,:\,}\tilde G^{(h)}+\frac1h\eta(h\tilde G^{(h)}){\,:\,}\tilde G^{(h)},
$$
where the last term on the right-hand side can be bounded by
$$
\frac1h|\eta(h\tilde G^{(h)}){\,:\,}\tilde G^{(h)}|\le \omega(h \|\tilde G^{(h)}\|_{L^\infty})|\tilde G^{(h)}|^2
$$
with $\omega(t)\to 0$, as $t\to 0^+$.
Together with (\ref{tildeG}), we obtain  for 
every $\varphi\in C^\infty([0,L])$
\begin{equation} \label{Taylortilde}
\int_\Omega \varphi \leb \tilde G^{(h)}{\,:\,}\tilde G^{(h)}\, dx
- 
\int_\Omega \varphi  \tilde E^{(h)} {\,:\,} \tilde G^{(h)}\, dx
\to 0.
\end{equation}
We now show that
\begin{equation} \label{tildeenergy}
\int_\Omega \varphi \tilde{E}^{(h)}{\,:\,} \tilde{G}^{(h)} \, dx -
\int_\Omega \varphi E^{(h)}{\,:\,} G^{(h)} \, dx \to 0.
\end{equation}
Together with the convergence of energy (\ref{en-conv}), the weak convergence
of $\tilde{G}^{(h)}$ to $G$ and (\ref{Taylortilde}) this
implies that
$$
\lim_{h\to 0} \int_{\Omega}\varphi\leb(\tilde G^{(h)}-G){\,:\,}(\tilde G^{(h)}-G)\, dx
= 0
$$
for every $\varphi\in C^\infty([0,L])$
with $\varphi(0) = 0$, and hence
\begin{equation}\label{conv-tildeG}
\sym (\tilde G^{(h)}-G) \to 0 \qquad \hbox{strongly in } 
L^2((a,L){\times}(-\tfrac12, \tfrac12);\mtwo)
\end{equation}
for all $a > 0$.
Using again a Taylor expansion we easily deduce that
\begin{equation} \label{conv-tildeE}
\tilde{E}^{(h)} \to E
 \qquad \hbox{strongly in } 
L^2((a,L){\times}(-\tfrac12, \tfrac12);\mtwo).
\end{equation}
To prove (\ref{tildeenergy}), we write the difference as 
\begin{equation}\label{diffe}
\int_{\Omega}\varphi  E^{(h)}{\,:\,} (\tilde G^{(h)}- G^{(h)})\, dx
+ \int_{\Omega} \varphi  ( \tilde E^{(h)} -E^{(h)} )
{\,:\,}\tilde G^{(h)}\, dx.
\end{equation}
The first  term can be controlled by
 the div-curl lemma; indeed, equalities (\ref{defGtilde}) and (\ref{form1}) yield
$$
R^{(h)}(\tilde G^{(h)}- G^{(h)})= \nabla_h(\tilde z^{(h)} -z^{(h)}),
$$
so that, by the Euler-Lagrange equation (\ref{EL2}), we have
$$
\begin{array}{c}
\displaystyle
\int_{\Omega} \varphi E^{(h)}
{\,:\,}(\tilde G^{(h)}- G^{(h)})\, dx
=  \int_{\Omega} \varphi R^{(h)}E^{(h)}
{\,:\,}\nabla_h(\tilde z^{(h)} - z^{(h)})\, dx
\smallskip
\\
\displaystyle
= h\int_{\Omega}\varphi\, g{\,\cdot\,} (\tilde z^{(h)} - z^{(h)})\, dx
-\int_{\Omega} \varphi'R^{(h)}E^{(h)}e_1{\,\cdot\,} (\tilde z^{(h)} - z^{(h)})\, dx.
\end{array}
$$
Since the sequence $\tilde z^{(h)} - z^{(h)}$ converges to $0$ strongly in $L^2(\Omega;\R^2)$ and
$R^{(h)}E^{(h)}$ is bounded in $L^2(\Omega;\mtwo)$, we conclude that
$$
\lim_{h\to 0} 
\int_{\Omega} \varphi E^{(h)}{\,:\,}(\tilde G^{(h)}- G^{(h)})\, dx
= 0.
$$
To estimate the second term in (\ref{diffe}),
%%%MGM:added reference
we use that
$\tilde E^{(h)}$ and $E^{(h)}$ are bounded in $L^2(\Omega;\mtwo)$. 
Thus using H\"older's inequality, we find
$$
\int_{\Omega} |\varphi (\tilde E^{(h)}- E^{(h)})
{\,:\,}\tilde G^{(h)}| \, dx \le \Big( \int_{A_h} |\tilde G^{(h)}|^2\, dx \Big)^{1/2}
\leq C \left[( \lambda_h^2 + h^{-1}) \leb^2(A_h)\right]^{1/2},
$$
and the right-hand side converges to zero in view of (\ref{tilde4}).
This concludes the proof of (\ref{tildeenergy}).

\medskip

\noindent {\em Step 3. Passage to the limit in the Euler-Lagrange
equation.}
\smallskip

\noindent
To pass to the limit in (\ref{EL5}), i.e., in the equation
\begin{equation}\label{EL5bis}
\int_0^L \Big( \hat E^{(h)}_{11}\phi' +(\theta^{(h)})'\hat E^{(h)}_{21}\phi
+\tfrac1h  \bar E^{(h)}_{21}\phi +\tfrac1h (\bar E^{(h)}_{12} - \bar E^{(h)}_{21}) \phi \Big)\, dx_1 =0,
\end{equation}
for  all $\phi\in C^\infty([0,L])$ which vanishes on an interval $(0,a)$, 
we first prove that
\begin{equation}\label{a1}
\lim_{h\to 0}\int_0^L (\theta^{(h)})'\hat E^{(h)}_{21}\phi \, dx_1 =0.
\end{equation}
Indeed, it follows from the definition of $\hat E^{(h)}$ that
\begin{eqnarray}
\lefteqn{ \int_0^L (\theta^{(h)})'\hat E^{(h)}_{21}\phi \, dx_1 =
\int_\Omega x_2(\theta^{(h)})' E^{(h)}_{21}\phi \, dx}
\nonumber
\\
& = & \int_\Omega x_2(\theta^{(h)})' \tilde E^{(h)}_{21}\phi \, dx 
+ \int_{A_h} x_2(\theta^{(h)})' (E^{(h)}_{21}- \tilde E^{(h)}_{21} )\phi \, dx.
\label{a02}
\end{eqnarray}
Using (\ref{conv-tildeE}) and the fact that $E_{21}=0$ a.e.\ in $\Omega$, 
we obtain that the first integral of the right-hand side converges to $0$. 
As for the second term, using H\"older's inequality, we obtain
$$
\int_{A_h} | x_2(\theta^{(h)})' (E^{(h)}_{21}- \tilde E^{(h)}_{21} )\phi | \, dx \le
C\Big( \int_{A_h} |(\theta^{(h)})' |^2\, dx  \Big)^{1/2}.
$$
Since $|(\theta^{(h)})'|\le Ch^{-1/2}$, and $h^{-1}\leb^2(A_h)\to 0$ by (\ref{tilde4}),
this implies that the second integral on the right-hand side of (\ref{a02})
also converges to $0$. This proves (\ref{a1}).

Using again a div-curl argument, we will  show that
\begin{equation}\label{b1}
\lim_{h\to 0}\int_0^L \tfrac1h (\bar E^{(h)}_{12} - \bar E^{(h)}_{21}) \phi \, dx_1= 0
\end{equation}
for every $\phi\in C^\infty([0,L])$ which vanishes on $(0,a)$.

Indeed, by (\ref{symm}) we have
$$
\tfrac1h \skw E^{(h)} = - \skw (E^{(h)}(G^{(h)})^T).
$$
Hence by (\ref{form1})
$$
\tfrac1h \skw E^{(h)} = - \skw ( E^{(h)} (\nabla_h z^{(h)})^T R^{(h)} ) +
x_2 (\theta^{(h)})' \skw ( E^{(h)} e_1\otimes e_1).
$$
We now use the fact that for $R \in SO(2)$ and $A \in \mtwo$ we have
$\skw A = R (\skw A) R^T = \skw (R A R^T)$ to deduce that
$$
\tfrac1h (E_{12}^{(h)} - E_{21}^{(h)}) = 
- (R^{(h)}E^{(h)})_{1j} (\nabla_h z^{(h)})_{2j}
+ (R^{(h)}E^{(h)})_{2j} (\nabla_h z^{(h)})_{1j}
- x_2 (\theta^{(h)})' E_{21}^{(h)}.
$$
Thus 
\begin{eqnarray*}
& & \int_0^L \tfrac1h (\bar{E}_{12}^{(h)} - \bar{E}_{21}^{(h)}) \phi \, dx_1 \\
&= & \int_\Omega 
 \left( -(R^{(h)}E^{(h)})_{1j} (\nabla_h z^{(h)})_{2j}
+ (R^{(h)}E^{(h)})_{2j} (\nabla_h z^{(h)})_{1j} \right) \phi \, dx \\
& & -\int_0^L  (\theta^{(h)})' \hat{E}_{21}^{(h)} \phi \, dx_1.
\end{eqnarray*}
The last term converges to zero by (\ref{a1}).
The Euler-Lagrange equation (\ref{EL2})
 and the strong convergence of $z^{(h)}$ imply that the remaining terms 
converge to
\begin{eqnarray*}
& \displaystyle  \int_\Omega  (R E)_{11} z_2 \phi' -  
(RE)_{21}  z_1 \phi'
\, dx \\
&\displaystyle  = 
\int_0^L  (R \bar{E})_{11} z_2 \phi' -
(R \bar{E})_{21}  z_1  \phi'\, dx_1  = 0,
\end{eqnarray*}
where we have used the fact that $z$ is independent of $x_2$
and $\bar{E} = 0$.
This proves (\ref{b1}) and together with (\ref{a1}) this shows that we can pass to the 
limit in (\ref{EL5bis}). Thus, we obtain again (\ref{ELfin})
and the proof can be concluded as before.
\end{proof}

\section{Auxiliary results} 
\label{Sec:aux}

In this section, we collect and prove some auxiliary results
which were needed 
in the proof of Theorem~\ref{thm:1.1}. 

We begin with an approximation result for deformations having elastic energy of order $h^2$ by means
of smooth rotations. This is the point where the rigidity theorem
 by Friesecke, James, and 
M\"uller \cite{FJM02} plays a crucial role
(note that in two dimensions the proof of this result can be
streamlined using complex variables). 

\begin{proposition}\label{rotation}
Let $(u^{(h)})\subset W^{1,2}(\Omega;\R^2)$ be a sequence such that
$$
F^{(h)}(u^{(h)}) := \int_\Omega \dist^2(\nabla_h u^{(h)}, SO(2))\, dx\le Ch^2,
$$
for every $h>0$. 
Then there exists an associated sequence $(R^{(h)})\subset C^\infty((0,L);\mtwo)$ such that 
\begin{eqnarray}
& \displaystyle 
R^{(h)} (x_1)\in SO(2) \quad \hbox{for every }x_1\in(0,L), \label{rot1}
\\
& \displaystyle 
\|\nabla_h u^{(h)} - R^{(h)}\|_{L^2} \le Ch, \label{rot2}
\\
& \displaystyle 
\|(R^{(h)})'\|_{L^2}+
h\, \| (R^{(h)})''\|_{L^2}\le  C  \label{rot3} 
\end{eqnarray}
for every $h>0$. If, in addition, 
$u^{(h)}(0,x_2) = (0, h x_2)$, then
\begin{equation} \label{bcrigid}
| R^{(h)}(0) - Id | \leq C \sqrt{h}.
\end{equation}
\end{proposition}

\begin{proof}
The argument follows closely the proof of Theorem~2.1 in \cite{Mor-Mue}. We include the details
for the convenience of the reader. 

For every $h>0$ let $k_h$ be an integer such that
$h\leq \frac{L}{k_h} <2h$. For every $a\in [0,L)\cap \frac{L}{k_h}\N$, we define 
$$
I_{a, h}:=(a, a+\tfrac{L}{k_h}).
$$
We apply the rigidity estimate \cite[Theorem~3.1]{FJM02} to the functions $v^{(h)}(z_1,z_2):=u^{(h)}(z_1, z_2/h)$ restricted to the set $(a,a+2h){\times}(-\frac h2,\frac h2)$ when $a< L- \frac{L}{k_h}$ and
 restricted to the set $(L-2h,L){\times}(-\frac h2,\frac h2)$, otherwise. 
Thus, we obtain
piecewise constant maps $Q^{(h)}: [0,L]\to SO(2)$ such that
\begin{equation}\label{01bis}
\int_{I_{a,h}{\times}(-\frac12,\frac12)}|\nabla_h u^{(h)} - Q^{(h)} |^2dx \leq C 
\int_{(a,a+2h){\times}(-\frac12, \frac12)} \dist^2(\nabla_h u^{(h)}, SO(2))\, dx.
\end{equation}
When $a=L-\frac{L}{k_h}$, replace the interval $(a,a+2h)$ by $(L-2h,L)$ in the second integral
above. We point out that the constant $C$ above is independent of~$h$.
Summing over $a$, we obtain
\begin{equation}\label{01}
\int_{\Omega}|\nabla_h u^{(h)} - Q^{(h)} |^2dx \leq C F^{(h)}(u^{(h)}) \leq Ch^2.
\end{equation}

Let $a\in [0,L)\cap \frac{L}{k_h}\N$ be such that $(a,a+4h)\subset (0,L)$ and let
$b:=a+ \frac{L}{k_h}$. Then, using estimate (\ref{01bis}), its analog for the set $(a,a+4h){\times}(-\frac12, \frac12)$, and the fact that both intervals $I_{a,h}$, $I_{b,h}$ are contained in $(a,a+4h)$, we have
\begin{equation}\label{diffquot}
\frac{L}{k_h}|Q^{(h)}(a) - Q^{(h)}(b)|^2
\leq C \int_{(a,a+4h){\times}(-\frac12,\frac12)} \dist^2(\nabla_h u^{(h)}, SO(2))\, dx.
\end{equation}
In particular, we deduce that
\begin{equation}\label{diffquot1}
|Q^{(h)}(x_1+s) - Q^{(h)}(x_1)|^2 \leq Ch^{-1}F^{(h)}(u^{(h)}) \leq Ch
\end{equation}
for every $x_1\in (h,L-h)$ and every $|s|\leq h$.
 If we extend $Q^{(h)}$ by $Q^{(h)}(0)$
for $x_1 < 0$ and by $Q^{(h)}(L)$ for $x_1 > L$, 
estimate (\ref{diffquot1}) holds for all $x_1 \in \R$.

Iterative application of inequality (\ref{diffquot}) provides 
a difference quotient estimate for~$Q^{(h)}$. More precisely,
\begin{equation}\label{11}
\int_\R  |Q^{(h)}(x_1+s)- Q^{(h)}(x_1)|^2dx_1 \leq 
C h^{-2} ( |s| + h )^2 F^{(h)}(u^{(h)})\leq C ( |s| + h )^2.
\end{equation}
Let $\eta\in C^{\infty}_0(0,1)$ be such that $\eta\geq 0$ and $\int_0^1\eta(s)\,ds=1$.
We set $\eta_h(s):=\frac{1}{h}\eta(\frac{s}{h})$ and define
$$
\tilde Q^{(h)}(x_1):= \int_{0}^h \eta_h(s) Q^{(h)}(x_1-s)\, ds, \qquad x_1\in[0,L].
$$
Using estimate (\ref{11}), we easily see that
\begin{equation}\label{001}
\begin{array}{c}
\|\tilde Q^{(h)} - Q^{(h)}\|_{L^2}\leq Ch,
\smallskip
\\
\|(\tilde Q^{(h)})'\|_{L^2}\le C, \qquad
h\, \| (\tilde Q^{(h)})''\|_{L^2}\le  C 
\end{array}
\end{equation}
for every $h>0$.

Let $\pi:U\to SO(2)$ be a smooth projection from a neighbourhood $U$ of $SO(2)$ onto $SO(2)$.
Since from (\ref{diffquot1}) and Jensen's  inequality we have that
\begin{equation}\label{inftyest}
\|\tilde Q^{(h)} - Q^{(h)} \|^2_{L^{\infty}}\leq Ch,
\end{equation}
the functions $\tilde Q^{(h)}$ take values in $U$ for $h$ small enough. Therefore, we can define $R^{(h)}:=\pi(\tilde Q^{(h)})$. Properties (\ref{rot1})--(\ref{rot3}) follow immediately from (\ref{01}) and~(\ref{001}).

%%% NEW Feb 3, 2006
To establish (\ref{bcrigid}), we start from 
the following trace inequality
\begin{equation}
\int_{(-\frac12, \frac12)} |v(0,z_2) - \bar{v}|^2 \,  dz_2
\leq
C \int_{(0,l){\times}(-\frac12,\frac12)}  |\nabla v|^2 \, dz, 
\end{equation}
which holds uniformly for $ 1 \leq l \leq 2$, with
$ \bar{v} = \int  v(0,z_2) \, dz_2$. We apply this estimate
with $v(z) = h^{-1} u^{(h)}(h z_1, z_2) - \bar{Q} z$, where
$\bar{Q} = Q^{(h)}(0)$. 
In combination with (\ref{01bis}) at $a= 0$ and the boundary condition for $u^{(h)}$,
this yields
\begin{equation}
 \int_{(-\frac12, \frac12)} \Big(|\bar{Q}_{12} x_2|^2 +
 |x_2 -   \bar{Q}_{22}  x_2 |^2\Big) \, dx_2
\leq C h.
\end{equation}
Since $\bar{Q} \in SO(2)$, this implies that $|\bar{Q} -Id | \leq C \sqrt{h}$.
In view of (\ref{inftyest}), this yields (\ref{bcrigid}).
\end{proof}

\begin{proposition} \label{Ehweak}
Assume that the energy density $W$ is differentiable and 
its derivative $DW$ is Lipschitz continuous. Assume moreover
that $DW$ is differentiable at the identity.
Suppose that 
$$ G^{(h)} \wto G \quad \mbox{weakly in } L^2(\Omega;\mtwo)
$$
and define the rescaled stresses as in  (\ref{Ehdef}) by
$$ E^{(h)} := \frac1h DW(Id + h G^{(h)}). 
$$
Then
\begin{equation}
E^{(h)}\wto E:=\leb\, G \quad \mbox{weakly in }L^2(\Omega;\mtwo),
\end{equation}
where $\leb:=D^2W(Id)$.
\end{proposition}

\begin{proof}
Since $(E^{(h)})$ is bounded in $L^2(\Omega;\mtwo)$, it is enough to show that the limit of each weakly convergent subsequence of $(E^{(h)})$ coincides with $\leb G$. Therefore, let $E^{(h_k)}$ converge weakly in $L^2(\Omega;\mtwo)$ to some~$\tilde E$.
 
A Taylor expansion of $DW$ around the identity yields
\begin{equation}\label{Taylor}
\frac1h DW(I+hA) = \leb A+ \frac1h \eta (hA),
\end{equation}
where $|\eta(F)|/|F| \to 0$ as $|F| \to 0$. 
Set $\omega(t):=\sup\{ |\eta(F)|/|F|:\ |F|\leq t\}$ for every $t>0$. 
Then $\omega(t)\to 0$, as $t\to 0^+$.
Applying formula (\ref{Taylor}) with $A$ replaced by $G^{(h_k)}$ and $h$ replaced by $h_k$, we find that
\begin{equation}\label{Taylor-k}
E^{(h_k)}=\frac{1}{h_k} DW(Id +{h_k} G^{(h_k)})= \leb G^{(h_k)} +\frac{1}{h_k}\, \eta(h_kG^{(h_k)})
\end{equation}
for every $k$.
Now let $M_k:=\{x\in\Omega: |G^{(h_k)}(x)|\le h_k^{-1/2}\}$,
and let $\chi_k$ be its characteristic function. Since $\chi_k\to 1$ 
boundedly in measure and $G^{(h)}\wto G$ weakly in $L^2(\Omega;\mtwo)$,
we have
\begin{equation}\label{weaks}
\chi_k E^{(h_k)}\wto \tilde E \quad \hbox{and} \quad \chi_k G^{(h_k)}\wto G \quad \hbox{weakly in }L^2(\Omega;\mtwo).
\end{equation} 
Moreover, from (\ref{Taylor-k}) it follows that for every $k$
\begin{equation}\label{Taylor-chi}
\chi_k E^{(h_k)} =  \leb(\chi_kG^{(h_k)}) +\frac{1}{h_k}\chi_k\, \eta(h_kG^{(h_k)}).
\end{equation}
The first term on the right-hand side converges 
weakly in $L^2(\Omega;\mtwo)$ to $\leb G$ by (\ref{weaks}). 
As for the second term, we have
$$
\frac{1}{h_k}\chi_k\, | \eta(h_kG^{(h_k)})|\le \omega(\sqrt{h_k})|G^{(h_k)}|,
$$
where we used the fact that $h_k|G^{(h_k)}|\le \sqrt{h_k}$ on $M_k$. Since $\omega(\sqrt{h_k})$ converges to~$0$ and $|G^{(h_k)}|$ is bounded in $L^2(\Omega;\mtwo)$, we conclude that the second term on the right-hand side of (\ref{Taylor-chi}) converges to $0$ strongly in $L^2(\Omega;\mtwo)$.
Passing to the limit in (\ref{Taylor-chi}) and using (\ref{weaks}), we finally obtain $\tilde E=\leb G$, so that the proposition is proved.
\end{proof}

We conclude this section by proving 
the truncation lemma stated and used in Section~\ref{sec:trcp}.
We first consider a fixed domain in $\R^n$ and then extend  the
result to 
the thin domains $\Omega_h$ by successive reflection.

\begin{lemma} \label{goodtrunc}
Let $1 < p < \infty$, $n \geq 2$, $m \geq 1$, $n,m  \in \N$ 
and let $\Omega$ be
a bounded Lipschitz domain in $\R^n$ or $\Omega = \R^n$. Then
there exists a constant $C_2$ (depending on $\Omega$, $p$, $n$, and $m$)
such that for every $A > a > 0$ and every function 
$v: \Omega \to \R^m$ with $\nabla v \in 
L^p(\Omega; {\mathbb M}^{m \times n})$
there exist $\lambda \in [a,A]$ and a function $v^\lambda: \Omega
\to \R^m$ such that
\begin{eqnarray}
|\nabla v^\lambda| &\leq& \lambda,  \label{gtrunc1} \\
\lambda^p | \{ v^\lambda \neq v \} | & \leq&
\frac{C_2}{\ln (A/a)} \int_\Omega  |\nabla v|^p \, dx.
\label{gtrunc2}
\end{eqnarray}
\end{lemma}
\begin{remark}
Estimate (\ref{gtrunc2}) is useful, if $A$ is much larger than $a$.
If $A$ and $a$ are comparable one can use  the simpler estimate
$\lambda^p |\{ v^\lambda \neq v \} | \leq C_2 ||\nabla v||_p^p$ 
for all $\lambda > 0$,
which even holds for $p=1$, see, e.g., \cite{EG}.
\end{remark}
\begin{proof} This is a simple consequence of the 
standard results on the truncation of gradients, see, e.g.,
\cite{Liu, Zi89, EG}. We sketch the argument for the convenience
of the reader. For simplicity, we only consider the scalar case $m=1$,
the general case can easily be treated by using the maximal function of 
$|\nabla v|$ and extending each component of $v$ separately. 
We first observe that it suffices to prove that there exists constants
$C'_1$ and $C'_2$ such that one can always find $\lambda \in [a,A]$ with
the following properties
\begin{eqnarray}
|\nabla v^\lambda| &\leq & C'_1 \lambda, \label{ggtrunc1} \\
\lambda^p |\{ v^\lambda \neq v \} | &\leq & \frac{C'_2}{\ln (A/a)} 
\int_\Omega |\nabla v|^p \, dx. \label{ggtrunc2}
\end{eqnarray}
Indeed, if (\ref{ggtrunc1}) and (\ref{ggtrunc2}) hold for $w$ and $w^\lambda$
(in place of $v$ and $v^\lambda$)
take $w = C'_1 v$ and $v^\lambda := w^{\lambda}/C'_1$. 
Then we obtain (\ref{gtrunc1}) and (\ref{gtrunc2}) with $C_2 := C'_2 (C'_1)^p$.

We first consider the case $\Omega = \R^n$. Let $f$ denote the 
Hardy-Littlewood maximal function of $\nabla v$,
$$ f(x) := \sup_{R >0} \frac{1}{|B(x,R)|} \int_{B(x,R)} |\nabla v| \, dy
$$
and let
$$
E^t := \{ x \in \R^n : f(x) > t \}.
$$
By standard estimates for the maximal function, we have
$$
||f||_{L^p} \leq C_3 ||\nabla v||_{L^p},
$$
where $C_3$ depends only on $p$. From the definition of $E^t$ 
and Poincar\'e inequality, one can easily deduce that
(see, e.g.,  \cite{EG}, p. 253)
%%%MGM:the lipschitz constant is C_4 t
$$
|v(x) - v(y)| \leq C_4 t \,|x-y|, \quad \mbox{for a.e. } x,y \not\in E^t.
$$
Hence (after removal of a null set), $v$ is Lipschitz in the complement
of $E^t$ and hence has an extension $v^t$ with the same Lipschitz 
constant. Thus
\begin{eqnarray*}
|\nabla v^t| &\leq & C_4 t, \\
|\{ v \neq v^t \}| & \leq & |E^t|.
\end{eqnarray*}
{}From the definition of $E^t$, we have the trivial estimate
$t^p |E^t| \leq ||f||_p^p \leq C_3^p ||\nabla v||_p^p$.
To obtain the refined estimate (\ref{gtrunc2}),
%%%MGM:corrected reference
we use the relation
$$
\int_0^\infty p t^{p-1} |E^t| \, dt
= \int_{\R^n} f^p \, dx.
$$ 
Set $g(t) := t^p |E^t|$. Then
$$
 \int^{A}_a \frac{p}{t} \inf_{t\in [a,A]} g(t) \, dt 
\leq \int^{A}_a \frac{p}{t} g(t) \, dt
\leq \int_{\R^n} |f|^p \, dx 
\leq C_3^p \int_{\R^n} |\nabla v|^p \, dx.
$$
This yields (\ref{gtrunc2}) with $C_2(\R^n) = \frac1p C_3^p$.

Suppose that $\Omega$ is a bounded Lipschitz domain. We may assume
that $\int_\Omega v  = 0$ (otherwise we first define $w:= v - m$, where
$m$ is the average of $v$, apply the result for $w$ and finally set
$v^\lambda := w^\lambda + m$). Then Poincar\'e's inequality yields
$$
||v||_{W^{1,p}(\Omega)} \leq C_5 ||\nabla v||_{L^p(\Omega)}.
$$
Thus there exists an extension $\tilde{v}:\R^n \to \R^m$ with
(see, e.g., \cite{St70})
$$
||\tilde{v}||_{W^{1,p}(\R^n)} \leq C_6 ||v||_{W^{1,p}(\Omega)}
\leq C_7 ||\nabla v||_{L^p(\Omega)}.
$$
Now we can apply the previous reasoning to $\tilde v$ and we get
(\ref{gtrunc1}) and (\ref{gtrunc2}) with 
$C_2(\Omega) = C_2(\R^n) C_7^p$.
\end{proof}

\begin{proof}[Proof of Lemma~\ref{truncate}]
Let $h>0$ and let $u\in W^{1,2}(\Omega_h;\R^2)$. 
First of all, note that $u$ can be extended to the rectangle $\Omega=(0,L){\times}(-\frac12,\frac12)$ by successive 
reflection. 
By Lemma~\ref{goodtrunc} there exist $\lambda \in [a,A]$ 
and $w \in W^{1,\infty}(\Omega;\R^2)$ such that
\begin{equation}\label{wOm}
\|\nabla w\|_{L^\infty(\Omega)}\le \lambda
\end{equation}
and 
\begin{equation}\label{wOm2}
\lambda^2 \leb^2(\{x\in\Omega: \ u(x)\neq w(x)\}) \le 
\frac{C}{\ln(A/a)}
\int_{\Omega} |\nabla u|^2\, dx.
\end{equation}
Let $N_h$ be the largest integer such that $hN_h+h/2\le 1/2$. For $i\in \Z\cap [-N_h,N_h]$ let $S_{h,i}:=(0,L){\times}(ih-\frac h2,ih+ \frac h2)$ and let 
$$
R_h:= \Omega\setminus  \bigcup_{-N_h\le i\le N_h} S_{h,i}.
$$
Since
$$
\sum_{-N_h\le i \le N_h} \leb^2(\{u\neq w\}\cap S_{h,i})\le  \leb^2(\{u\neq w\}),
$$
there exists some index $i_0$ such that 
\begin{eqnarray}
%%%MGM:multiplied by lambda^2
\lambda^2 \leb^2(\{u\neq w\}\cap S_{h,i_0})
& \le & \frac{1}{2N_h+1} \lambda^2 \leb^2(\{u\neq w\}) \nonumber 
\\
& \le &
\frac{C}{2 N_h + 1} \frac{1}{\ln(A/a)} \int_\Omega |\nabla u|^2 \, dx.
\label{i0}
\end{eqnarray}
Let $v:\Omega_h\to\R^2$ be the function defined by
$$
v(x_1,x_2):= w(x_1, i_0h+(-1)^{i_0} x_2) \qquad \hbox{for every } (x_1,x_2)\in\Omega_h.
$$
It is clear that $v\in W^{1,\infty}(\Omega_h;\R^2)$ and that it satisfies (\ref{trunc1}) by (\ref{wOm}). Moreover, since $u$ has been extended to $\Omega$ by reflection, we have
\begin{equation} \label{i1}
\{x\in\Omega_h: \ u(x)\neq v(x)\} = \{ x\in S_{h,i_0}:\ u(x)\neq w(x)\}
\end{equation}
and
\begin{equation} \label{i2}
\int_\Omega |\nabla u|^2 \, dx \le
(2 N_h + 3) \int_{\Omega_h} |\nabla u|^2 \, dx.
\end{equation}
Now assertion (\ref{trunc2}) follows from (\ref{i0})--(\ref{i2}).
This finishes the proof of Lemma~\ref{truncate}.
\end{proof}

\section*{Acknowledgements}
We would like to thank R. Pakzad for many helpful comments.
M.G.M. and S.M. were supported by the Marie Curie
research training network MRTN-CT-2004-505226 (MULTIMAT). 
M.G.M. was also partially supported by MIUR project 
``Calculus of Variations" 2004.

\end{document}